# Weyl cycles on the blow-up of $\mathbb{P}^4$ at eight points

Maria Chiara Brambilla, Olivia Dumitrescu and Elisa Postinghel

**Abstract** We define the *Weyl cycles* on $X^n_s$, the blown up projective space $\mathbb{P}^n$ in $s$ points in general position. In particular, we focus on the Mori Dream spaces $X^3_7$ and $X^4_8$, where we classify all the Weyl cycles of codimension two. We further introduce the *Weyl expected dimension* for the space of the global sections of any effective divisor that generalizes the *linear expected dimension* of [2] and the *secant expected dimension* of [3].

## 1 Introduction

Let $X^n_s$ be the blown up projective space $\mathbb{P}^n$ in $s$ points in general position. When the number of points $s$ is small, the space $X^n_s$ has an interpretation as certain moduli space, see e.g. [1] and [4]. Mori Dream Spaces of the form $X^n_s$ were classified via the work of Mukai [16, 17] and techniques of birational geometry of moduli spaces. In previous work, in order to analyze properties of the pairs $(X^n_s, D)$ with $D$ a Cartier divisor, the authors of this article developed techniques of *polynomial interpolation theory* in [2] for $s = n + 2$ and in [3] for $s = n + 3$ respectively, via the study of the base loci.

Maria Chiara Brambilla
Università Politecnica delle Marche, via Brecce Bianche, I-60131 Ancona, Italy, e-mail: brambilla@dipmat.univpm.it

Olivia Dumitrescu
1. University of North Carolina at Chapel Hill, 340 Phillips Hall CB 3250 NC 27599-3250 and
2. Simion Stoilow Institute of Mathematics Romanian Academy 21 Calea Grivitei Street 010702 Bucharest Romania, e-mail: dolivia@unc.edu
The second author is supported by NSF grant DMS1802082,

Elisa Postinghel
Dipartimento di Matematica, Università degli Studi di Trento, via Sommarive 14 I-38123 Povo di Trento (TN), Italy e-mail: elisa.postinghel@unitn.it

In this paper we apply an analogous approach, based on interpolation theory, to define and study the subvarieties determining the birational geometry of the Mori Dream Spaces $X_7^3$ and $X_8^4$. We will use this study as an opportunity to reveal the geometry hidden by the Weyl group action on fixed linear cycles of $X_s^n$ and its consequences. For instance, we predict that for all Mori Dream Spaces, Weyl cycles determine the birational geometry of such spaces, cones of effective and movable divisors and their decomposition into nef chambers.

In this article we propose a definition of *Weyl cycles* on $X_s^n$ as follows (see Definition 1 for details). We call *Weyl divisor* any effective divisor $D$ in $\mathrm{Pic}(X_s^n)$ in the Weyl orbit of an exceptional divisor $E_i \in \mathrm{Pic}(X_s^n)$. We call a *Weyl cycle of codimension i* an element of the Chow group $A^i(X_s^n)$ that is an irreducible component of the intersection of Weyl divisors, which are pairwise orthogonal with respect to the Dolgachev-Mukai pairing on $\mathrm{Pic}(X_s^n)$.

For an arbitrary number of points, $s$, Weyl divisors are always generators of the Cox ring. The correspondence between $(-1)$-curves of $\mathbb{P}^2$ and Weyl curves in $X_s^2$ (i.e. Weyl divisors) was proved by Nagata [18], while giving a counterexample to the Hilbert 14-th problem. Moreover in the case of $\mathbb{P}^2$, Weyl curves have been widely investigated since they are involved in the well-known Segre-Harbourne-Gimigliano-Hirschowitz conjecture, see e.g. [5, 6] or, more widely, in the classification of algebraic surfaces (Castelnuovo's contraction theorem), the base of the minimal model program. The notion of *divisorial* $(-1)$-*classes* on $X_s^n$ was introduced by Laface and Ugaglia in [15] and recently studied by the second author and Priddis in [12].

For $X_s^3$, Laface and Ugaglia introduced the notion of elementary $(-1)$-curves and studied their properties in [14]. The case of Weyl cycles of $X_8^4$ has been studied in [10] with a different approach: indeed in such paper Weyl orbits in $X_8^4$ (as well as in $X_7^3$) of the proper transforms of linear cycles blown up along lines spanned by any two points are described. The classification of Weyl cycles obtained in [10] for the cases $X_7^3$, $X_8^4$, yields the same classification we determine here for Weyl curves in $X_7^3$ and for Weyl surfaces (see Equations (1)) in $X_8^4$. Therefore we conclude that Definition 1 and the definitions used in [10] are equivalent for cycles of codimension 2 in $X_7^3$ and $X_8^4$. The authors believe that these two definitions are related in general, namely for Weyl cycles in $X_s^n$, and they will study their connection in a forthcoming paper.

In this article we emphasize that basic methods of intersection theory, applied to pairs of orthogonal Weyl divisors, give an iterative method to compute the Weyl cycles of codimension 2 in $X_7^3$ and in $X_8^4$. Moreover, we show that every Weyl cycles is swept out by families of rational curves parametrized by Weyl cycles of larger codimension, see Proposition 3, Corollary 1 and Lemma 5.

Our main result, contained in Section 5.3, is a classification of all the Weyl surfaces of $X_8^4$. For each such surface, on the one hand we compute its multiplicity of containment in the base locus of any given divisor, on the other hand we compute its class in the Chow ring of $X_{8,(1)}^4$, the blow up of $X_8^4$ along the strict transforms of all lines through two base points and all rational normal quartic curves through seven base points; we believe that it is the only representative of its class. There are five such classes, up to index permutation, as listed in the following formula (see Section 5.2 for the precise notation):

$$
\begin{aligned}
S^1_{1,4,5} &: h - e_1 - e_4 - e_5 - \sum_{i,j \in \{1,4,5\}} (e_{ij} - f_{ij}) \\
S^3_{1,\widehat{8}} &: 3h - 3e_1 - \sum_{i=2}^{7} e_i - (e_{C_{\widehat{8}}} - f_{C_{\widehat{8}}}) - \sum_{i=2}^{7} (e_{1i} - f_{1i}) \\
S^6_{6,7,8} &: 6h - 3\sum_{i=1}^{5} e_i - \sum_{i=6}^{8} e_i - \sum_{i,j \in \{1,2,3,4,5\}, i \neq j} (e_{ij} - f_{ij}) - \sum_{k=6}^{8} (e_{C_{\widehat{k}}} - f_{C_{\widehat{k}}}) \\
S^{10}_{1,2} &: 10h - 6e_1 - 6e_2 - \sum_{i=3}^{8} 3e_i - 3(e_{12} - f_{12}) - \sum_{i=1}^{2} \sum_{j=3}^{8} (e_{ij} - f_{ij}) - \sum_{k=3}^{8} (e_{C_{\widehat{k}}} - f_{C_{\widehat{k}}}) \\
S^{15}_{8} &: 15h - \sum_{i=1}^{7} 6e_i - 3e_8 - \sum_{1 \leq i < j \leq 7} (e_{ij} - f_{ij}) - \sum_{i=1}^{7} (e_{C_{\widehat{i}}} - f_{C_{\widehat{i}}}) - 3(e_{C_{\widehat{8}}} - f_{C_{\widehat{8}}})
\end{aligned}
\tag{1}
$$

Recall that the birational geometry of $X^4_8$ has been investigated in [17] and [4]. Casagrande, Codogni and Fanelli studied in detail the relation between the geometry of $X^2_8$ and $X^4_8$ and in [4, Theorem 8.7] they described five types of surfaces in $X^4_8$ playing a special role in the Mori program. We emphasize that this list agrees with our classification of Weyl surfaces, (1).

Finally, we propose here a notion of expected dimension for a linear system which takes into account the contribution to the speciality given by the special cycles contained in the base locus. In Definition 2, we introduce the *Weyl expected dimension* for $X^3_7$ and $X^4_8$, which extends the analogous definitions of *linear expected dimension* of [2] and *secant expected dimension* of [3]. We prove that any effective divisor $D$ in $X^3_7$ satisfies $h^0(X^3_7, O_{X^3_7}(D)) = \text{wdim}(D)$, see Theorem 1, and we conjecture that the same holds in $X^4_8$.

The paper is organized as follows. In Section 2, we introduce the notation, recall basic facts on the blown up of $\mathbb{P}^n$ in $s$ general points, $X^n_s$, and on the action of standard Cremona transformations on $\text{Pic}(X^n_s)$. In Section 3 we introduce the definition of Weyl cycles and we give some general result on Weyl curves in $X^n_s$. Section 4 is devoted to the preliminary case of $X^3_7$, where we classify Weyl divisors and Weyl curves and we describe their geometry. Section 5 concerns the case of $X^4_8$. The main result, i.e. the classification of the Weyl surfaces is contained in Section 5.3. In Section 5.4, we give the classification of Weyl divisors and their geometrical description. The last Section 6 is devoted to the dimensionality problem.

**Acknowledgements:** The first and third authors are members of GNSAGA-INDAM. The second author is supported by the NSF grant DMS - 1802082. The third author was partially supported by the EPSRC grant EP/S004130/1.

We wish to thank Cinzia Casagrande for many useful discussions.

## 2 Preliminaries

We denote by $X_s^n$ the blown up of $\mathbb{P}^n$ at $s$ general points $\mathcal{I} = \{p_1, \ldots, p_s\}$. The Picard group of $X_s^n$ is $\text{Pic}(X_s^n) = \langle H, E_1, \ldots, E_s \rangle$, where $H$ is a general hyperplane class, and the $E_i$'s are the exceptional divisors of the $p_i$'s. For any subset $J \subseteq \{1, \ldots, s\}$ of cardinality $\leq n$, we denote by $L_J$ the class, in the Chow ring of $X_s^n$, of the strict transform of the linear cycle spanned by $J$. If $|J| = n$, then $L_I = H - \sum_{i \in J} E_i \in \text{Pic}(X_s^n)$ is the class of a fixed hyperplane.

The *Dolgachev-Mukai pairing* on $\text{Pic}(X_s^n)$ is the bilinear form defined as follows (cf. [16]):
$$\langle H, H \rangle = n - 1, \quad \langle H, E_i \rangle = 0, \quad \langle E_i, E_j \rangle = -\delta_{i,j}.$$

The *standard Cremona transformation based on the coordinate points* on $\mathbb{P}^n$ is the birational transformation defined by the following rational map:
$$\text{Cr} : (x_0 : \cdots : x_n) \to (x_0^{-1} : \cdots : x_n^{-1}),$$

see e.g. [9, 12] for more details. Given any subset $I \subseteq \{1, \ldots, s\}$ of cardinality $n+1$, we denote by $\text{Cr}_I$ and call *standard Cremona transformation* the map obtained by precomposing $\text{Cr}$ with a projective transformation which takes the points indexed by $I$ to the coordinate points of $\mathbb{P}^n$. A standard Cremona transformation induces an automorphism of $\text{Pic}(X_s^n)$, denoted again by $\text{Cr}_I$ by abuse of notation, by sending a divisor
$$D = dH - \sum m_i E_i \tag{2}$$
to
$$\text{Cr}_I(D) = (d - c)H - \sum_{i \in I}(m_i - c)E_i - \sum_{j \notin I}^{s} m_j E_j, \tag{3}$$

where $c := \sum_{i \in I} m_i - (n-1)d$. The canonical divisor $-(n+1)H + (n-1)\sum_{i=1}^{s} E_i$ is invariant under such an automorphism. The *Weyl group* $W_{n,s}$ acting on $\text{Pic}(X_s^n)$ is the group generated by standard Cremona transformations, see [9]. We say that a divisor (2) is Cremona reduced if $c \leq 0$ for any $I$ of cardinality $n+1$.

In [12, Theorem 3.2] the authors observed that the intersection pairing between divisors is preserved under Cremona transformation.

**Lemma 1** *Let $D, F$ be two divisors and let $\omega \in W_{n,s}$ be an element of the Weyl group. Then $\langle \omega(D), \omega(F) \rangle = \langle D, F \rangle$.*

Here we point out that the scheme-theoretic intersection of two divisors is in general not preserved under Cremona transformation. Let $D, F$ be two divisors and let $\omega = \text{Cr}_I$ be a standard Cremona transformation. Then
$$\omega(D \cap F) \cup \Lambda \subseteq \omega(D) \cap \omega(F)$$
where $\Lambda$ is a union of linear cycles of the indeterminacy locus of $\omega$. The following lemma provides an explicit recipe for $\Lambda$.

**Lemma 2** *Let $I \subseteq \{1, \ldots, s\}$ have cardinality $n + 1$, and let $I = I_1 \cup I_2$, with $|I_1| = m + 1$ and $|I_2| = n - m$. Let $D = dH - \sum m_i E_i$ be a divisor in $X_s^n$. If $(n - m - 1)d - \sum_{i \in I_2} m_i = a \geq 1$, then the $m$-plane $L_{I_1}$ is contained in $\mathrm{Cr}_I(D)$ exactly $a$ times.*

*Proof* Set $c = \sum_{i \in I} m_i - (n-1)d$. By [2] and [11, Proposition 4.2], we can compute the multiplicity of containment of the $m$-plane $L_{I_1}$ in $\mathrm{Cr}(D)$:

$$\sum_{i \in I_1}(m_i - c) - m(d - c) = \sum_{i \in I_1} m_i - md - c = (n - m - 1)d - \sum_{i \in I_2} m_i = a,$$

concluding the proof. □

## 3 Weyl cycles in $\mathbb{P}^n$ blown up at $s$ points

In [12, Definition 4.1] a smooth divisor $D$ in $\mathrm{Pic}(X_s^n)$ is called $(-1)$-*class* (or $(-1)$-*divisorial cycle*) if $D$ is effective, integral and it satisfies $\langle D, D \rangle = -1$ and $\langle D, -K_{X_s^n} \rangle = n - 1$. In [12, Theorem 0.5], it is proved that $D$ is a $(-1)$-class if and only if it is in the Weyl orbit of some exceptional divisor $E_i$. Notice that if $i \in I$, then $\mathrm{Cr}_I(E_i) = L_{I \setminus \{i\}}$ is a hyperplane through $n$ base points.

Here we generalize the definition of $(-1)$-classes to cycles of higher codimension in $X_s^n$, as follows. We will say that two divisors $D$ and $F$ are *orthogonal* if $\langle D, F \rangle = 0$.

**Definition 1** We say that an effective divisor $D \in \mathrm{Pic}(X_s^n)$ is a *Weyl divisor* if it belongs to the Weyl orbit of an exceptional divisor $E_i$. A non-trivial effective cycle $C \in A^i(X_s^n)$ is a *Weyl cycle of codimension $i$* if it is an irreducible component of the intersection of pairwise orthogonal Weyl divisors.

*Remark 1* Let $s \geq n + 1$ and $1 \leq m \leq n - 1$. Any $m$-plane $L$ spanned by $m + 1$ points is a Weyl cycle. Indeed, it is easy to check that $L$ is the intersection of $r = n - m$ pairwise orthogonal hyperplanes spanned by $n$ points. By Lemma 1, any effective cycle $C$ contained in the Weyl orbit of a $m$-plane $L$ spanned by $m + 1$ base points is a Weyl cycle. In particular the Weyl planes and Weyl lines studied in [10] are always Weyl cycles, according to Definition 1.

### 3.1 Weyl curves

We collect here some results on Weyl cycles of codimension $n - 1$ in $X_s^n$, which we call *Weyl curves*. The following examples show explicitly that the strict transforms of lines through two points and of the rational normal curves of degree $n$ through $n + 3$ points are Weyl curves in $X_s^n$, according to Definition 1.

*Example 1* Let $L = L_{12}$ be the line through $p_1$ and $p_2$, then $L = D_1 \cap \cdots \cap D_{n-1}$ where $D_i = L_{I_i}$ and $I_i = \{1, 2, \ldots, n + 1\} \setminus \{i + 2\}$ for any $1 \leq i \leq n - 1$.

*Example 2* For any $i = 1, \ldots, n-1$, consider the pairwise orthogonal Weyl divisors $D_i = 2H - 2E_1 - \ldots - 2E_{n-1} - E_n - E_{n+1} - E_{n+2} - E_{n+3} + E_i$. Once can easily check that $D_1 \cap \cdots \cap D_{n-1}$ is the union of $L_{1\ldots n-1}$ and the rational normal curve of degree $n$ through $n+3$ points.

We recall that the Chow group of algebraic curves $A^{n-1}(X_s^n)$ is generated by $h^1, e_i^1$, the classes of a general line in $X_s^n$ and of a general line on the exceptional divisor $E_i$, respectively. When $n = 4$, the following formula describes the action on curves of the standard Cremona transformation $\mathrm{Cr}_J$, based on the set $J$ [10] for $n = 4$ if $C = \delta h^1 - \sum_{i=1}^{s} \mu_i e_i^1$, then

$$\mathrm{Cr}_J(C) = (4\delta - 3\sum_{j \in J} \mu_j)h^1 - \sum_{j \in J}(\delta - \sum_{i \in J \setminus \{j\}} \mu_i)e_j^1 - \sum_{j \notin J} \mu_j e_j^1. \qquad (4)$$

Case $n = 3$ is described in [13]. The general formula given in equation 4 will be proved in an upcoming paper.

*Remark 2* Given a divisor $D$ in $X_s^n$ and a line $L_{ij} = h^1 - e_i^1 - e_j^1$, then the multiplicity of containment of the line $L_{ij}$ in the base locus of $D$ is exactly $\max\{0, -D \cdot L_{ij}\}$, where $\cdot$ denotes the intersection product in the Chow ring of $X_s^n$ (cf. [11, Proposition 4.2]). If $n = 3, 4$, the same holds for any curve $C$ in the Weyl orbit of the line $L_{ij}$, thanks to formulas (4).

## 4 $\mathbb{P}^3$ blown up in 7 points

In this section we consider Weyl cycles of $X_7^3$, the blow up of $\mathbb{P}^3$ at 7 points in general position. Recall that $X_7^3$ is a Mori Dream Space and that the cone of effective divisors is generated by the divisors of anticanonical degree $\frac{1}{2}\langle D, -K_{X_7^3} \rangle = 1$. These are exactly the Weyl divisors and they fit in five different types, modulo index permutation.

**Proposition 1** *The Weyl divisors in $X_7^3$ are, modulo index permutation:*

*(1) $E_i$ (exceptional divisor);*
*(2) $H - E_1 - E_2 - E_3$ (planes through three points);*
*(3) $2H - 2E_1 - E_2 - E_3 - E_4 - E_5 - E_6$ (pointed cone over the twisted cubic);*
*(4) $3H - 2(E_1 + E_2 + E_3 + E_4) - E_5 - E_6 - E_7$ (Cayley nodal cubic);*
*(5) $4H - 3E_1 - 2(E_1 + E_2 + E_3 + E_4 + E_5 + E_6 + E_7)$.*

*Proof* It is easy to compute the Weyl orbit of a plane through 3 points, by applying formula (3). □

**Proposition 2** *The Weyl curves in $X_7^3$ are the fixed lines $h^1 - e_i^1 - e_j^1$ and the fixed twisted cubics $3h^1 - \sum_{i=1}^{7} e_i^1 + e_j^1$.*

*Proof* For every pair of orthogonal Weyl divisors as in Proposition 2, one can check that the intersection is the union of fixed lines and twisted cubics. □

From the previous result we can conclude that our Definition 1 of Weyl curves in $X_7^3$ is equivalent to the definition of Weyl line of [10].

In the following result we describe the intrinsic geometry of the Weyl divisors of $X_7^3$, showing that they are covered by pencils of rational curves parametrized by a Weyl curve.

**Proposition 3** *Let $D$ be Weyl divisor on $X_7^3$. If $C \subset D$ is a Weyl curve, then there is a pencil of rational curves $\{C_q : q \in C\}$ with $C_q \cdot D = 0$ sweeping out $D$.*

*Proof* We will consider the divisors (2)-(5) from Proposition 2. It is easy to check what Weyl curves are contained in $D$, using Remark 2. For each such containment $C \subseteq D$, we will find a suitable pencil of curves, parametrised by $C$, sweeping out $D$.

(2) Let us consider the fixed hyperplane $D = H - E_1 - E_2 - E_3$ and the Weyl line $L_{12} \subset D$. Such plane is swept out by the pencil of lines through $p_3$ and with a point $q \in L_{12}$: $\{C_3^1(q) : q \in L_{12}\}$. Since the cycle class of $C_3^1(q)$ is $h^1 - e_3^1$, then we obtain $C_3^1(q) \cdot D = 0$.

(3) The quadric surface $D = 2H - 2E_1 - E_2 - E_3 - E_4 - E_5 - E_6$ contains the fixed twisted cubic $C_{1,\ldots,6}^3 = 3h^1 - \sum_{i=1}^6 e_i^1$. Since it is the strict transform of a pointed cone, it is swept out by the pencil of lines $\{C_1^1(q) : q \in C_{1,\ldots,6}\}$. Notice also that $D$ can be obtained from $H - E_1 - E_2 - E_3$ through the transformation $\mathrm{Cr}_{1,4,5,6}$ (cf. (3)). The latter preserves the line $L_{12}$ and, for every $q \in L_{12}$, it sends the line $C_1^1(q)$ to the cubic curve $C_{1,3,4,5,6}^3(q)$, see formula (4). Therefore we see that $D$ is also swept out by the pencil $\{C_{1,3,4,5,6}^3(q) : q \in L_{12}\}$.

(4) This surface is obtained from (3) via the standard Cremona transformation $\mathrm{Cr}_{2347}$. The image of the first pencil sweeping out (3) is the pencil of cubics $\{C_{1,\ldots,4,7}^3(q) : q \in C_{1,\ldots,6}\}$ and it sweeps out (4). The images of the second pencil sweeping out (3) is $\{C_{3,4}^5(q) : q \in L_{12}\}$, where $C_{3,4}^5(q)$ is a quintic curve with cycle class $5h - \sum_{i=1}^7 e_i - e_3 - e_4$ and passing through $q \in L_{12}$: this pencil sweeps out (4).

(5) This surface is obtained from (4) via $\mathrm{Cr}_{1,5,6,7}$, On the one hand we obtain that (5) is swept out by the pencil of quintics $\{C_{1,7}^5(q) : q \in C_{1,\ldots,6}\}$. On the other hand the surface is covered by the pencil of septic curves $C_2^7(q)$ with class $7h - 2\sum_{i=1}^7 e_i + e_2$ passing through $q \in L_{12}$: $\{C_2^7(q) : q \in L_{12}\}$.

In each of the above, since the intersection product $\cdot$ is preserved under Cremona transformation, we can conclude.

## 5 $\mathbb{P}^4$ blown up in 8 points

### 5.1 Curves in $X_8^4$

**Notation 1** We consider the following classes of moving curves in $A^3(X_8^4)$, each obtained from the previous via a standard Cremona transformation (see formula (4) and a permutation of indices. They each live in a 4-dimensional family.

- $h^1 - e_i^1$, for any $i \in \{1, \ldots, 8\}$,
- $4h^1 - \sum_{i \in J} e_i^1$, for any $J \subset \{1, \ldots, 8\}$ with $|J| = 6$,
- $7h^1 - \sum_{i \in J} 2e_i^1 - \sum_{i \notin J} e_i^1$, for any $J$ with $|J| = 3$,
- $10h^1 - e_{i_1} - 3e_{i_2} - \sum_{i \neq i_1, i_2} 2e_i^1$, for any $i_1 \neq i_2$, $i_1, i_2 \in \{1, \ldots, 8\}$,
- $13h^1 - \sum_{i \in J} 2e_i^1 - \sum_{i \notin J} 3e_i^1$, for any $J$ with $|J| = 3$.
- $16h^1 - \sum_{i \in J} 4e_i^1 - \sum_{i \notin J} 3e_i^1$, for any $J$ with $|J| = 2$.

The families of curves in Notation 1 correspond to facets of the effective cone of divisors on $X_8^4$, see [4]. Here we include a proof via our geometrical approach.

**Proposition 4** Let $D = dH - \sum m_i E_i$ be a divisor in $X_8^4$. If $D$ is effective, then we have:

- $m_i \leq d$, for every $i \in \{1, \ldots, 8\}$,
- $\sum_{i \in J} m_i - 4d \leq 0$, for any $J \subset \{1, \ldots, 8\}$ with $|J| = 6$,
- $\sum_{i \in J} 2m_i + \sum_{i \notin J} m_i - 7d \leq 0$, for any $J$ with $|J| = 3$,
- $m_{i_1} + 3m_{i_2} + \sum_{i \neq i_1, i_2} 2m_i - 10d \leq 0$, for any $i_1 \neq i_2$, $i_1, i_2 \in \{1, \ldots, 8\}$,
- $\sum_{i \in J} 2m_i + \sum_{i \notin J} 3m_i - 13d \leq 0$, for any $J$ with $|J| = 3$.
- $\sum_{i \in J} 4m_i + \sum_{i \notin J} 3m_i - 16d \leq 0$, for any $J$ with $|J| = 2$.

The first two inequalities were also proved in [2, Lemma 2.2].

**Proof** Notice that each 4-dimensional family of Notation 1 covers $X_8^4 \setminus \bigcup_i E_i$, indeed for each general point in $X_8^4 \setminus \bigcup_i E_i$ we find one curve of the family that passes through it. Now, if $D \cdot (h^1 - e_i^1) = d - m_i < 0$, then $D$ contains each line in the family in its base locus, but this contradicts the assumption that $D$ is effective. This proves the first inequality. The remaining inequalities are proved similarly. □

## 5.2 Further blow up of $\mathbb{P}^4$.

For any $1 \leq i \leq 8$, we denote by $C_{\hat{i}}$ the rational normal quartic curve passing through seven base points and skipping the $i$th point. Consider now

$$X_{8,(1)}^4 \xrightarrow{p} X_8^4,$$

the blow up of $X_8^4$ along the 28 lines $L_{ij}$ and the 8 curves $C_{\hat{i}}$. The Picard group of $X_{8,(1)}^4$ is $\text{Pic}(X_{8,(1)}^4) = < H, E_i, E_{ij}, E_{C_{\hat{i}}} >$, where, abusing notation, we denote again by $E_i$ the pull-back $p^*(E_i)$ and by $H$ the pull-back $p^*(H)$, while $E_{ij}$ and $E_{C_{\hat{i}}}$ are the exceptional divisors of the curves. Notice that $E_i$ is a blown up $\mathbb{P}^3$ in 14 general points, $X_{14}^3$, while $E_{ij} \cong \mathbb{P}^1 \times \mathbb{P}^2$ and $E_{C_{\hat{i}}} \cong C_{\hat{i}} \times \mathbb{P}^2$.

For any $D \in \text{Pic}(X_8^4)$, of the form $D = dH - \sum_{i=1}^{8} m_i E_i$, the strict transform $\widetilde{D}$ of $D$ under $p$ satisfies

$$\widetilde{D} := D - \sum k_{ij} E_{ij} - \sum k_{C_{\hat{i}}} E_{C_{\hat{i}}}. \tag{5}$$

where $k_{ij}$ and $k_{C_{\hat{i}}}$ are defined in Remark 2.

Let us consider now the Chow group of 2-cycles of $X^4_{8,(1)}$:

$$A^2(X^4_{8,(1)}) = \langle h, e_i, e_{ij}, f_{ij}, e_{C_{\hat{i}}}, f_{C_{\hat{i}}} \rangle.$$

where $h$ is the pullback of a general plane of $\mathbb{P}^4$, $e_i$ is the pull-back of a general plane contained in $E_i$, $f_{ij} \cong \mathbb{P}^2$ is the fiber over a point of the line and $e_{ij} \cong \mathbb{P}^1 \times \mathbb{P}^1$ is the transverse direction, $f_{C_{\hat{i}}}$ is the fiber over a point of the curve $C_{\hat{i}}$ and $e_{C_{\hat{i}}}$ is the transverse direction. In the Chow ring $A^*(X^4_{8,(1)})$ we have the following relations:

$$H^2 = h, \quad E_i^2 = -e_i, \quad HE_i = 0, \quad E_iE_j = 0 \tag{6}$$

$$HE_{ij} = E_iE_{ij} = f_{ij}, \quad E_iE_{jk} = 0 \tag{7}$$

$$E_{ij}^2 = -e_{ij} - f_{ij}, \quad E_{ij}E_{ik} = 0, \quad E_{ij}E_{kl} = 0 \tag{8}$$

$$H^4 = h^2 = 1, \quad E_i^4 = e_i^2 = -1, \quad f_{ij}e_{ij} = e_{ij}^2 = -1 \tag{9}$$

$$f_{ij}^2 = e_if_{ij} = 0, \quad he_i = hf_{ij} = 0, \quad e_ie_{ij} = he_{ij} = 0. \tag{10}$$

### 5.3 Classification of the Weyl surfaces

The section contains one of the main results of this paper. We construct five Weyl surfaces in $X^4_8$ and we prove that they are the only such cycles, modulo index permutation. For any surface, we also give its exact multiplicity of containment in a given divisor and its class in the Chow ring of $X^4_{8,(1)}$.

**Proposition 5** *Let $S^1 = S^1_{1,4,5}$ be the plane $L_{145}$ through three points in $\mathbb{P}^4$.*

- *Given an effective divisor $D = dH - \sum m_i E_i$ in $X^4_8$, let*

$$k_{S^1}(D) = \max\{0, m_1 + m_4 + m_5 - 2d\}.$$

*Then the surface $S^1$ is contained in the base locus of $D$ exactly $k_{S^1}(D)$ times.*

- *The class of the strict transform $\widetilde{S^1}$ of $S^1$ in the Chow group $A^2(X^4_{8,(1)})$ is*

$$h - e_1 - e_4 - e_5 - \sum_{i,j \in \{1,4,5\}} (e_{ij} - f_{ij}).$$

*Proof* The first part of the statement follows from [2] and [11, Proposition 4.2].

Consider the fixed hyperplanes $D_0 := H - E_1 - E_3 - E_4 - E_5$ and $F_0 := H - E_1 - E_2 - E_4 - E_5$. Let $\widetilde{D_0}$ and $\widetilde{F_0}$ be their strict transforms on $X^4_{8,(1)}$, see (5). Clearly we have $S^1 = D_0 \cap F_0$, and $\widetilde{S^1} = \widetilde{D_0} \cap \widetilde{F_0}$. By using relations (6), (7), (8), we compute $\widetilde{D_0} \cap \widetilde{F_0} = h - e_1 - e_4 - e_5 - \sum_{i,j \in \{1,4,5\}} (e_{ij} - f_{ij})$. □

Using Lemma 2 we obtain the following.

**Lemma 3** *Given a subset $I = \{i_1, \ldots, i_5\} \subseteq \{1, \ldots, 8\}$ and a divisor $D = dH - \sum m_i E_i$ in $X_8^4$. If*
$$d - m_{i_1} - m_{i_2} = a \geq 1,$$
*then the 2-plane $L_{i_3 i_4 i_5}$ is contained in $\mathrm{Cr}_I(D)$ exactly $a$ times.*

**Lemma 4** *Let $I = \{i_1, i_2, i_3, i_4, i_5\}$ and $J = \{i_1, i_2, i_6\}$ two subset of $\{1, \ldots, 8\}$, such that $|I \cap J| = 2$. If $\mathrm{Cr}_I$ is the standard Cremona transformation based on I, then the plane $L_J$ is $\mathrm{Cr}_I$-invariant, that is $\mathrm{Cr}_I(L_J) = L_J$.*

*Proof* Consider the hyperplanes $D = L_{i_1 i_2 i_3 i_6}$ and $F = L_{i_1 i_2 i_4 i_6}$. We have $D \cap F = L_{i_1 i_2 i_6} = L_J$. Clearly $\mathrm{Cr}_I(D) = D$ and $\mathrm{Cr}_I(F) = F$ and hence also
$$\mathrm{Cr}_I(L_J) = \mathrm{Cr}(D \cap F) \subseteq \mathrm{Cr}_I(D) \cap \mathrm{Cr}_I(F) = D \cap F = L_J.$$

**Proposition 6** *Let $J := \{1, 2, 3, 6, 7\}$ and consider the Cremona transformation $\mathrm{Cr}_J$. Let $S^1 = L_{145}$. Then $S^3 := \mathrm{Cr}_J(S^1)$ is the strict transform of cubic pointed cone over the rational normal curve $C_{\widetilde{8}}$ and the point $p_1$.*

- *Given a divisor $D = dH - \sum m_i E_i$, let*
$$k_{S^3}(D) = \max\{0, 2m_1 + m_2 + m_3 + m_4 + m_5 + m_6 + m_7 - 5d\}.$$
*Then the surface $S_3$ is contained in the base locus of $D$ exactly $k_{S^3}(D)$ times.*

- *The class of the strict transform $\widetilde{S^3}$ of $S^3$ in the Chow group $A^2(X_{8,(1)}^4)$ is*

$$3h - 3e_1 - \sum_{i=2}^{7} e_i - (e_{C_{\widetilde{8}}} - f_{C_{\widetilde{8}}}) - \sum_{i=2}^{7}(e_{1i} - f_{1i})$$

*Proof* The plane $S^1 = L_{145}$ is swept out by the pencil of lines $\{C^1(q) : q \in L_{14}\}$, where the cycle class of $C^1(q)$ is $h^1 - e_5^1$ and it passes through the point $q \in L_{14}$. Using formulas (4) and the same idea as in the proof of Proposition 3, we compute the images of the line $L_{14} = h^1 - e_1^1 - e_4^1$ and of the pencil of lines $\{C^1(q) : q \in L_{14}\}$ of class $h^1 - e_5^1$ via the transformation $\mathrm{Cr}_J$. We have $\mathrm{Cr}_J(L_{14}) = L_{14}$ and $\mathrm{Cr}_J(C^1(q)) = C^4(q)$ where $C^4(q)$ is rational curve with class $4h^1 - e_1^1 - e_2^1 - e_3^1 - e_5^1 - e_6^1 - e_7^1$ and passing through $q$. Thus we get that the surface $S_3$ is swept out by the pencil $\{C^4(q) : q \in L_{14}\}$. Therefore $D$ contains any curve $C^4(q)$, and hence $S_3$, in its base locus at least $\max\{0, m_1 + m_2 + m_3 + m_5 + m_6 + m_7 + \max\{0, m_1 + m_4 - d\} - 4d\}$ times. Notice that we have $m_1 + m_2 + m_3 + m_5 + m_6 + m_7 - 4d \leq 0$, since $D$ is effective, by Proposition 4. Hence the claim follows.

Now we prove the second statement. Given $D_0$ and $F_0$ defined in the previous proposition, recall that $D_0 \cap F_0 = S^1$. We consider now their images $D_1 = \mathrm{Cr}_J(D_0)$ and $F_1 = \mathrm{Cr}_J(F_0)$:

$$D_1 = 2H - 2E_1 - E_2 - 2E_3 - E_4 - E_5 - E_6 - E_7$$
$$F_1 = 2H - 2E_1 - 2E_2 - E_3 - E_4 - E_5 - E_6 - E_7.$$

Clearly $S^3 \subseteq D_1 \cap F_1$. By Proposition 5 we easily see that the only plane contained in $D_1 \cap F_1$ is $L_{123}$. Moreover it is easy to check that the intersection $D_1 \cap F_1$ does not intersect the indeterminacy locus of the Cremona transformation $\text{Cr}_J$ in any other 2-dimensional component. Hence $D_1 \cap F_1$ is the union of the plane $L_{123}$ and an irreducible cubic surface with one triple point in $p_1$ and 6 simple points. We conclude that $S^3$ is exactly such cubic surface.

We now shall describe the class of $S^3$ in $A^2(X^4_{8,(1)})$. Let $\widetilde{D_1}$ and $\widetilde{F_1}$ be the corresponding strict transforms under the blow up of lines and rational normal curves in $X^4_{8,(1)}$. By (5), we have

$$\widetilde{D_1} = 2H - 2E_1 - E_2 - 2E_3 - \sum_{i=4}^{7} E_i - 2E_{13} - \sum_{i \in \{1,3\}, k \in \{2,4,5,6,7\}} E_{ik} - E_{C_{\hat{8}}}$$

$$\widetilde{F_1} = 2H - 2E_1 - 2E_2 - \sum_{i=3}^{7} E_i - 2E_{12} - \sum_{1 \le i \le 2, 3 \le k \le 7} E_{ik} - E_{C_{\hat{8}}}.$$

By using relations (6), (7), (8), we compute the intersection:

$$\widetilde{D_1} \cap \widetilde{F_1} = (h - e_1 - e_2 - e_3 - \sum_{i,j \in \{1,2,3\}} (e_{ij} - f_{ij})) + (3h - 3e_1 - \sum_{i=2}^{7} e_i - (e_{C_{\hat{8}}} - f_{C_{\hat{8}}}) - \sum_{i=2}^{7} (e_{1i} - f_{1i})).$$

Finally by Proposition 5, we can conclude. □

We will denote by $S^3_{i,\hat{j}}$ the cubic surface with a triple point at $p_i$ and multiplicity zero at $p_j$.

**Proposition 7** *Let* $J := \{2, 3, 4, 5, 8\}$ *and consider the Cremona transformation* $\text{Cr}_J$. *Then* $S^6 := \text{Cr}_J(S^3)$ *is a surface of degree 6 with five triple points.*

- *Given an effective divisor* $D = dH - \sum m_i E_i$, *let*

$$k_{S_6}(D) = \max\{0, 2(m_1 + m_2 + m_3 + m_4 + m_5) + m_6 + m_7 + m_8 - 8d\}.$$

  *Then the surface* $S_6$ *is contained in the base locus of* $D$ *exactly* $k_{S_6}(D)$ *times.*
- *The class in* $A^2(X^4_{8,(1)})$ *of strict transform* $\widetilde{S^6}$ *of* $S^6$ *in* $X^4_{8,(1)}$ *is*

$$6h - 3\sum_{i=1}^{5} e_i - \sum_{i=6}^{8} e_i - \sum_{i,j \in \{1,2,3,4,5\}, i \ne j} (e_{ij} - f_{ij}) - \sum_{k=6}^{8} (e_{C_{\hat{k}}} - f_{C_{\hat{k}}}).$$

*Proof* We know from the previous proposition that the surface $S^3$ is swept out by a pencil of rational normal quartic curves $\{C^4(q) : q \in L_{14}\}$. By (4), we obtain that the image of the pencil is $\{C^7(q) : q \in L_{14}\}$, where $C^7(q)$ is a rational septic curve with class $7h^1 - \sum_{i=1}^{8} e_i^1 - e_2^1 - e_3^1 - e_5^1$ and passing through $q \in L_{14}$. Since the surface $S^6$ is swept out by this, we can say that $D$ contains $S^6$ in its base locus at least $\max\{0, m_1 + 2m_2 + 2m_3 + m_4 + 2m_5 + m_6 + m_7 + m_8 + \max\{0, m_1 + m_4 - d\} - 7d\}$ times. Since $D$ is effective, by Proposition 4 we have $m_1 + 2m_2 + 2m_3 + m_4 + 2m_5 + m_6 + m_7 + m_8 - 7d \le 0$, hence the claim follows.

Now we prove the second statement. Given $D_1$ and $F_1$ defined in the previous proposition, recall that $D_1 \cap F_1 = L_{123} \cup S^3$. We consider now $D_2 = \mathrm{Cr}_J(D_1)$ and $F_2 = \mathrm{Cr}_J(F_1)$ to be their image under the Cremona transformation and we get:

$$D_2 = 3Hx - 2E_1 - 2E_2 - 3E_3 - 2E_4 - 2E_5 - E_6 - E_7 - E_8$$
$$F_2 = 3H - 2E_1 - 3E_2 - 2E_3 - 2E_4 - 2E_5 - E_6 - E_7 - E_8.$$

We now analyse the intersection $D_2 \cap F_2$. Note that $\mathrm{Cr}_J(L_{123}) = L_{123}$, by Lemma 4. By Proposition 5 we see that the only planes contained in $D_2 \cap F_2$ are $L_{123}, L_{234}, L_{235}$. Finally we check that the intersection of $D_2 \cap F_2$ with the indeterminacy locus of $\mathrm{Cr}_J$ does not contain any other 2-dimensional component, besides the planes $L_{234}$ and $L_{235}$. Hence the intersection $D_2 \cap F_2$ splits into four components: the three planes $L_{123}, L_{234}, L_{235}$ and a sextic surface with five triple points at $p_1$ and three simple points. Hence we conclude that $S^6$ is exactly the sextic irreducible surface.

We now describe the class of $S^6$ in $X_{8,(1)}^4$. Let $\widetilde{D_2}$ and $\widetilde{F_2}$ be the corresponding strict transforms under the blow up of lines and rational normal curves in $X_{8,(1)}^4$, see (5). We have

$$\widetilde{D_2} = 3H - \sum_{i\in\{1,2,4,5\}} 2E_i - 3E_3 - \sum_{i=6}^{8} E_i - \sum_{i\in\{1,2,4,5\}} 2E_{3i} - \sum_{i=6}^{8} E_{3i} - \sum_{i,j\in\{1,2,4,5\}, i\neq j} E_{ij} - \sum_{i=6}^{8} E_{C_{\hat{i}}}$$

$$\widetilde{F_2} = 3H - \sum_{i\in\{1,3,4,5\}} 2E_i - 3E_2 - \sum_{i=6}^{8} E_i - \sum_{i\in\{1,3,4,5\}} 2E_{2i} - \sum_{i=6}^{8} E_{2i} - \sum_{i,j\in\{1,3,4,5\}, i\neq j} E_{ij} - \sum_{i=6}^{8} E_{C_{\hat{i}}}$$

Computing their complete intersection, we have:

$$\widetilde{D_2} \cap \widetilde{F_2} = (h - e_1 - e_2 - e_3 - \sum_{i,j\in\{1,2,3\}, i\neq j}(e_{ij} - f_{ij})) +$$

$$+ (h - e_2 - e_3 - e_4 - \sum_{i,j\in\{2,3,4\}, i\neq j}(e_{ij} - f_{ij})) + (h - e_2 - e_3 - e_5 - \sum_{i,j\in\{2,3,5\}}(e_{ij} - f_{ij})) +$$

$$(6h - 3\sum_{i=1}^{5} e_i - \sum_{i=6}^{8} e_i - \sum_{i,j\in\{1,2,3,4,5\}, i\neq j}(e_{ij} - f_{ij}) - \sum_{k=6}^{8}(e_{C_{\hat{k}}} - f_{C_{\hat{k}}}))$$

where we use relations (6), (7), (8), and we conclude. □

We will denote by $S^6_{i,j,k}$ the sextic surface with five triple points at $\{p_h\}$ for $h \neq i, j, k$.

**Proposition 8** *Let* $J := \{1, 2, 6, 7, 8\}$ *and consider the Cremona transformation* $\mathrm{Cr}_{lJ}$. *Then* $S^{10} := \mathrm{Cr}_J(S^6)$ *is a surface of degree* 10 *with two sextuple points and six triple points.*

- *Given an effective divisor $D = dH - \sum m_i E_i$, let*

$$k_{S^{10}}(D) = \max\{0, 3(m_1 + m_2) + 2(m_3 + m_4 + m_5 + m_6 + m_7 + m_8) - 11d\}.$$

*Then the surface $S^{10}$ is contained in the base locus of $D$ exactly $k_{S^{10}}(D)$ times.*
- *The class of the strict transform $\widetilde{S^{10}}$ of $S^{10}$ in $A^2(X^4_{8,(1)})$ is*

$$10h - 6e_1 - 6e_2 - \sum_{i=3}^{8} 3e_i - 3(e_{12} - f_{12}) - \sum_{i=1}^{2}\sum_{j=3}^{8}(e_{ij} - f_{ij}) - \sum_{k=3}^{8}(e_{C_{\widehat{k}}} - f_{C_{\widehat{k}}}).$$

*Proof* We know from the previous proposition that the surface $S^6$ is swept out by the pencil of rational septic curves $\{C^7(q) : q \in L_{14}\}$. By (4), we obtain that the image of the pencil is $\{C^{10}(q) : q \in L_{14}\}$, where $C^{10}(q)$ is a rational curve with class $10h^1 - 2e_1^1 - 3e_2^1 - 2e_3^1 - e_4^1 - 2e_5^1 - 2e_6^1 - 2e_7^1 - 2e_8^1$ and passing through $q \in L_{14}$. Since the surface $S^{10}$ is swept out by this pencil, we can say that $D$ contains $S^{10}$ in its base locus at least $\max\{0, 2m_1 + 3m_2 + 2m_3 + m_4 + 2m_5 + 2m_6 + 2m_7 + 2m_8 + \max\{0, m_1 + m_4 - d\} - 10d\}$ times. Since $2m_1 + 3m_2 + 2m_3 + m_4 + 2m_5 + 2m_6 + 2m_7 + 2m_8 - 10d \leq 0$, by Proposition 4 then we the claim follows.

Now we prove the second statement. Given $D_2$ and $F_2$ defined in the previous proposition, recall that $D_2 \cap F_2 = S^6 \cup L_{123} \cup L_{234} \cup L_{235}$. We consider now $D_3 = \text{Cr}_J(D_2)$ and $F_3 = \text{Cr}_J(F_2)$ to be their image under the Cremona transformation and we get

$$D_3 = 5H - 4E_1 - 4E_2 - 3E_3 - 2E_4 - 2E_5 - 3E_6 - 3E_7 - 3E_8$$
$$F_3 = 4H - 3E_1 - 4E_2 - 2E_3 - 2E_4 - 2E_5 - 2E_6 - 2E_7 - 2E_8.$$

It is easy to check, by applying the previous propositions, that the intersection $D_3 \cap F_3$ contains the planes $L_{123}, L_{126}, L_{127}, L_{128}$ and the cubic surfaces $S^3_{2,\widehat{4}}$ and $S^3_{2,\widehat{5}}$. Notice that $\text{Cr}_J(L_{123}) = L_{123}$ by Lemma 4, and $\text{Cr}_J(L_{234}) = S^3_{2,\widehat{5}}$, $\text{Cr}_J(L_{235}) = S^3_{2,\widehat{4}}$, by Proposition 6. By computing the intersection of $D_3 \cap F_3$ with the indeterminacy locus of $\text{Cr}_J$ we see that there are no other 2-dimensional components, besides the planes $L_{126}, L_{127}, L_{128}$. Hence we conclude that $S^{10}$ is an irreducible surface with degree 10 and two sixtuple points at $p_1$ and $p_2$ and 6 triple points.

Finally we describe the class of $S^{10}$ in $X^4_{8,(1)}$. Let $\widetilde{D_3}$ and $\widetilde{F_3}$ be the corresponding strict transforms under the blow up of lines and rational normal curves in $X^4_{8,(1)}$, see (5). Computing their complete intersection, as in the previous case we get our claim. □

**Proposition 9** *Let $J := \{3, 4, 5, 6, 7\}$ and consider the Cremona transformation $\text{Cr}_J$. Then $S^{15} := \text{Cr}_J(S^{10})$ is a surface of degree 15 with one triple point and seven sextuple points.*

- *Given an effective divisor $D = dH - \sum m_i E_i$, let*

$$k_{S^{15}}(D) = \max\{0, 3(m_1 + m_2 + m_3 + m_4 + m_5 + m_6 + m_7) - 2m_8 - 14d\}.$$

*Then the surface $S^{15}$ is contained in the base locus of D exactly $k_{S^{15}}(D)$ times.*
- *The class of the strict transform $\widetilde{S^{15}}$ of $S^{15}$ in $A^2(X^4_{8,(1)})$ is*

$$15h - \sum_{i=1}^{7} 6e_i - 3e_8 - \sum_{1 \le i < j \le 7} (e_{ij} - f_{ij}) - \sum_{i=1}^{7}(e_{C_{\widehat{i}}} - f_{C_{\widehat{i}}}) - 3(e_{C_{\widehat{8}}} - f_{C_{\widehat{8}}}).$$

*Proof* We know from the previous proposition that the surface $S^{10}$ is swept out by the pencil of rational septic curves $\{C^{10}(q) : q \in L_{14}\}$. By (4), we obtain that the image of the pencil is $\{C^{13}(q) : q \in L_{14}\}$, where $C^{13}(q)$ is a rational curve with class $13h^1 - 2e_1^1 - 3e_2^1 - 3e_3^1 - 2e_4^1 - 3e_5^1 - 3e_6^1 - 3e_7^1 - 2e_8^1$ and passing through $q \in L_{14}$. Since the surface $S^{15}$ is swept out by this pencil, we can say that $D$ contains $S^{15}$ in its base locus at least $\max\{0, 2m_1 + 3m_2 + 3m_3 + 2m_4 + 3m_5 + 3m_6 + 3m_7 + 2m_8 + \max\{0, m_1 + m_4 - d\} - 13d\}$ times. Since $2m_1 + 3m_2 + 3m_3 + 2m_4 + 3m_5 + 3m_6 + 3m_7 + 2m_8 - 13d \le 0$ by Proposition 4, then we the claim follows.

Now we prove the second statement. Given $D_3$ and $F_3$ defined in the previous proposition, recall that

$$D_3 \cap F_3 = S_{10} \cup L_{123} \cup L_{126} \cup L_{127} \cup L_{128} \cup S^3_{2,\widehat{4}} \cup S^3_{2,\widehat{5}}.$$

We consider now $D_4 = \mathrm{Cr}_J(D_3)$ and $F_4 = \mathrm{Cr}_J(F_3)$ to be their image under the Cremona transformation.

$$D_4 := 7H - 4E_1 - 4E_2 - 5E_3 - 4E_4 - 4E_5 - 5E_6 - 5E_7 - 3E_8$$
$$F_4 := 6H - 3E_1 - 4E_2 - 4E_3 - 4E_4 - 4E_5 - 4E_6 - 4E_7 - 2E_8.$$

Now the intersection $D_4 \cap F_4$ contains $S^3_{3,\widehat{8}} = \mathrm{Cr}_J(L_{123})$ (by Proposition 6), $S^3_{6,\widehat{8}} = \mathrm{Cr}_J(L_{126})$ (by Proposition 6), $S^6_{148} = \mathrm{Cr}_J(S^3_{2,\widehat{4}})$ (by Proposition 7), $S^6_{158} = \mathrm{Cr}_J(S^3_{2,\widehat{5}})$ (by Proposition 7). Moreover we have the components: $S^3_{7,\widehat{8}}$, $S^6_{128}$, and it can be easily proved that $S^3_{7,\widehat{8}} = \mathrm{Cr}_J(L_{127})$ and $S^6_{128} = \mathrm{Cr}_J(L_{128})$. Finally we check that the intersection of $D_4 \cap F_4$ with the indeterminacy locus of $\mathrm{Cr}_J$ does not contain any 2-dimensional component. Hence we conclude that $S^{15}$ is an irreducible surface of degree 15 and with a triple point at $p_8$ and seven sextuple points.

Finally, as in the previous case, we compute the complete intersection of the strict transforms $\widetilde{D_4}$ and $\widetilde{F_4}$, and we get our statement. □

*Remark 3* We point out that the five Weyl surfaces described above correspond to the same list computed by [4, Theorem 8.7].

*Remark 4* Notice that the class of a Weyl surface must be computed in the Chow ring of $X^4_{8,(1)}$ in order to uniquely identify the surface. Indeed, for instance, the class of $S^3$ in the Chow ring of $X^4_8$ is $3h - 3e_1 - \sum_{i=2}^{7} e_i$, but so is the class of the union of the three planes $L_{123}$, $L_{145}$ and $L_{167}$. This was predicted in [7, Theorem 4.4] where the authors proved that the cone of effective 2-cycles of $X^4_8$ is linearly generated, namely each effective cycle can be written as a sum of linear cycles. However, the three planes do not contain the rational normal curve, whereas $S^3$ does. From this observation it is clear that the cone of effective cycles of codimension 2 of $X^4_{8,(1)}$ will not be linearly generated.

*Remark 5* Notice that, in Propositions 5,6,7,8,9 we used a specific sequence of Cremona transformations to obtain each Weyl surface of $X^4_8$ from the previous. This choice is clearly not unique, in fact there are multiple paths going from one Weyl surface to another. Similarly, for each Weyl surface $S$ we found a suitable pencil of curves over a Weyl curve $C \subseteq S$ that covers it. This description is also not unique, in particular for every Weyl curve $C \subseteq S$, we can find one such pencil.

**Proposition 10** *The five surfaces $S^1, S^3, S^6, S^{10}, S^{15}$ are the only Weyl surfaces in $X^4_8$.*

*Proof* The statement can be proved by direct inspection. In Proposition 11 below we classify all the Weyl divisors in $X^4_8$. Then we consider all the possible intersection of two orthogonal Weyl divisors and, by using Propositions 5, 6, 7, 8, 9, and computing degrees and multiplicities, we have checked that all the irreducible components of the intersections are surfaces of type $S^1, S^3, S^6, S^{10}, S^{15}$. □

By the previous proposition we conclude that any Weyl surface of $X^4_8$ is contained in the orbit of a plane through 3 points. Hence our Definition 1 of Weyl surface in this case coincide with the definition Weyl plane given in [10].

From the proofs of Propositions 5, 6, 7, 8, 9, we get the following consequence.

**Corollary 1** *Every Weyl surface on $X^4_8$ is swept out by a pencil of rational curves curve $\{C(q) : q \in C\}$ over a Weyl curve $C$.*

### 5.4 Weyl divisors.

Recall that $X^4_8$ is a Mori Dream Space and in particular the cone of effective divisors is finitely generated by the divisors of anticanonical degree $\frac{1}{3}\langle D, -K_{X^4_8}\rangle = 1$. A simple application of formula (3) gives the following classification of all the Weyl divisors in $X^4_8$; they are exactly the generators of the effective cone, see also [17].

**Proposition 11** *The Weyl divisors in $X^4_8$ are, modulo permutation of indices:*

*(1) $E_i$, (the exceptional divisor)*
*(2) $H - \sum_{i=1}^{4} E_i$, (hyperplane through four points);*
*(3) $2H - 2E_1 - 2E_1 - \sum_{i=3}^{7} E_i$, (quadric cone, join of a rational normal quartic and a line);*
*(4) $3H - \sum_{i=1}^{7} 2E_i$, (the 2-secant variety to a rational normal quartic);*
*(5) $3H - 3E_1 - \sum_{i=2}^{5} 2E_i - \sum_{i=6}^{8} E_i$, (cone on the Cayley surface of $\mathbb{P}^3$);*
*(6) $4H - \sum_{i=1}^{4} 3E_i - \sum_{i=5}^{7} 2E_i - E_8$, with $|J| = 4$ and $j \notin J$;*

*(7)* $4H - 4E_1 - 3E_2 - \sum_{i=3}^{8} 2E_i$, *(cone on a quartic surface of $\mathbb{P}^3$)*;
*(8)* $5H - 4E_1 - 4E_2 - \sum_{i=3}^{6} 3E_i - 2E_7 - 2E_8$;
*(9)* $6H - 5E_1 - \sum_{i=2}^{4} 4E_i - \sum_{i=5}^{8} 3E_i$;
*(10)* $6H - \sum_{i=1}^{6} 4E_i - 3E_7 - 2E_8$;
*(11)* $7H - \sum_{i=1}^{3} 5E_i - \sum_{i=4}^{7} 4E_i - 3E_8$;
*(12)* $7H - 6E_1 - \sum_{i=2}^{8} 4E_i$;
*(13)* $8H - 6E_1 - \sum_{i=2}^{6} 5E_i - 4E_7 - 4E_8$;
*(14)* $9H - \sum_{i=1}^{4} 6E_i - \sum_{i=5}^{8} 5E_i$;
*(15)* $10H - 7E_1 - \sum_{i=2}^{8} 6E_i$.

We conclude this section with the following geometrical descriptions of the Weyl divisors on $X_8^4$. As pencils of curves with cycle class as in Notation 1 sweep out Weyl surfaces of $X_8^4$, nets of such curves sweep out Weyl divisors.

**Lemma 5** *Let $D$ be a Weyl divisor on $X_8^4$ containing a Weyl surface $S$. Then there is a net of curves $\{C(q) : q \in S\}$ with $C(q) \cdot D = 0$ sweeping out $D$.*

***Proof*** Notice that every divisor (2)-(15) of Proposition 5.4 satisfies the hypotheses. For one such divisor, let $S \subset D$. By Propositions 5-6-7-8-9, we can find a sequence of standard Cremona transformations such that the image of $S$ is a plane $S_1$. Applying the same sequence of transformations to $D_1$, we obtain a Weyl divisor containing such plane. Modulo reordering the points, the possible outputs for the image of $D$ are the divisors (2), (3), (5), (7), (8), (9), (11) of Proposition 11. For each such output, we shall exhibit a sequence of Cremona transformations that preserve the plane $S_1$ and takes $D_1$ to a hyperplane containing $S_1$. Without loss of generality, we will assume that $S_1$ is the class of the plane passing through the first three points. In the following tables, for every (i), on the left hand side we will describe the class of the Weyl divisor and on the right hand side the class of the curve $C(q)$ of the net:

```
(11)  7|5 5 5 3 4 4 4 4    19|4 4 4 3 4 4 4 4
 (9)  6|5 4 4 3 4 3 3 3    16|4 3 3 3 2 3 3 3
 (8)  5|4 4 3 2 3 3 3 2    13|3 3 2 2 3 3 3 2
 (7)  4|4 3 2 2 2 2 2 2    10|3 2 1 2 2 2 2 2
 (5)  3|3 2 2 2 2 1 1 1     7|2 1 1 2 2 1 1 1
 (3)  2|2 2 1 1 1 1 1 0     4|1 1 0 1 1 1 1 0
 (2)  1|1 1 1 1 0 0 0 0     1|0 0 0 1 0 0 0 0
```

## 6 Weyl expected dimension

Let $n = 3, 4$. For $r \in \{1, 2, 3\}$, let $L_{I(r)}$ be a linear cycle of dimension $r$ spanned by $r + 1$ base points. Recall that $W_{n,n+4}$ denotes the Weyl group of $X_{n+4}^n$. Consider the following set of Weyl $r$-cycles: $W_n(r) := \{w(L_{I(r)}) : w \in W_{n,n+4}\}$, and let $k_A(D)$ denote the multiplicity of containment of the $r$-cycle $A$ in the base locus of the divisor $D$.

By Remark 2 we know that for any Weyl curve $A \in W_n(1)$, then $k_C = \max\{0, -D \cdot A\}$. For every Weyl divisor $A \in W_n(n-1)$ (i.e. those listed in Propositions 2 and 11), we have that $k_A = -\max\{0, \langle D, A \rangle\}$, see [2, Proposition 2.3] and [11, Proposition 4.2] for details. Finally, for $n = 4$, by the results of Section 5.3 we know that $W_4(2)$ is the set of the Weyl surfaces (i.e. those listed in equation (1)) and the multiplicity of containement $k_A(D)$ of any Weyl surface $A \in W_4(2)$ in the base locus of an effective divisor $D$ is computed in Propositions 5, 6, 7, 8, 9.

We introduce now the notion of *Weyl expected dimension*.

**Definition 2** Let $n = 3, 4$ and $D$ be an effective divisor on $X = X^n_{n+4}$. We say that $D$ has *Weyl expected dimension* wdim($D$), where

$$\text{wdim}(D) := \chi(X, \mathcal{O}_X(D)) + \sum_{r=1}^{n-1} \sum_{A \in W_n(r)} (-1)^{r+1} \binom{n + k_A(D) - r - 1}{n}.$$

We now show that the Weyl expected dimension in invariant under the action of the Weyl group.

**Proposition 12** *Let $n = 3, 4$ and $D$ an effective divisor on $X^n_{n+4}$. The Weyl dimension of $D$ is preserved under standard Cremona transformations.*

*Proof* Let $D = dH - \sum_{i=1}^{n+4} m_i E_i$. We need to prove that wdim($D$) = wdim($\text{Cr}_I(D)$) for $\text{Cr}_I$ a standard Cremona transformation. Let $D' = dH - \sum_{i \in I} m_i E_i$ be the divisor obtained from $D'$ by forgetting 3 points. From [2, Corollary 4.8, Theorem 5.3] we have that wdim($D'$) = wdim($\text{Cr}_I(D')$), where the formula wdim($D'$) only takes into account the Weyl cycles of $D$ based exclusively at the points parametrized by $I$ that are therefore fixed linear subspaces through base points.

We claim that, for all the remaining Weyl cycles $A$ of $D$, interpolating at least a point away from the indeterminacy locus and for which $k_A(D) \geq 1$, we have $k_A(D) = k_{\text{Cr}(A)}(\text{Cr}(D))$. If $A$ is a curve, the claim is true because $k_A(D) = -A \cdot D = -\text{Cr}(A) \cdot \text{Cr}(D)$. If $A$ is a divisor, the claim is true because $k_A(D) = -\langle A, D \rangle = -\langle \text{Cr}(A), \text{Cr}(D) \rangle$. It only remains to show the claim for $A = S$ a surface of $X^4_8$. It follows from the proofs of Propositions 5, 6, 7, 8, 9 and Remark 5 that for a Weyl curve $C \subseteq S$ such that $S$ is swept out by a pencil $\{C(q) : q \in C\}$, then $k_S(D) = -C(q) \cdot D + k_C(D)$. Since $D$ is effective, then $C(q) \cdot D \geq 0$ by Proposition 4, so $k_C(D) = -C \cdot D \geq 1$. Since $k_S(D) = -C(q) \cdot D - C \cdot D = -\text{Cr}(C_q) \cdot \text{Cr}(D) - \text{Cr}(C) \cdot \text{Cr}(D) = -\text{Cr}(C(q)) \cdot D + k_{\text{Cr}(C)}(\text{Cr}(D))$ and $\text{Cr}(D)$ is swept out by $\{\text{Cr}(C(q)) : q \in \text{Cr}(C)\}$, we conclude. $\square$

This yields an explicit formula for the dimension of any linear system in $X^3_7$.

**Theorem 1** *For any effective divisor $D \in Pic(X^3_7)$, we have*

$$\text{h}^0(X^3_7, \mathcal{O}_{X^3_7}(D)) = wdim(D).$$

**Proof** For the sake of simplicity, we will abbreviate $h^0(X_7^3, O_{X_7^3}(D))$ with $h^0(D)$. Consider a sequence of standard Cremona transformations which takes $D$ to a Cremona reduced divisor $D'$: it is well-known that $h^0(D) = h^0(D')$. By the previous proposition we have that $wdim(D) = wdim(D')$. Since $D'$ is Cremona reduced, by [8, Theorem 5.3] we know that $D'$ is linearly non-special, i.e. its dimension equals its linear expected dimension introduced in [2]: $h^0(D') = ldim(D') = wdim(D')$ where the last equality is easy to check for Cremona reduced divisors in $X_7^3$. Hence we conclude that $h^0(D) = wdim(D)$. □

For the case of $X_8^4$, we propose the following conjecture.

**Conjecture 1** *For any effective divisor $D \in Pic(X_8^4)$, we have*

$$h^0(X_8^4, O_{X_8^4}(D)) == wdim(D).$$

We claim that the notion of Weyl dimension extends both that of *linear expected dimension* of [2] and that of *secant expected dimension* of [3].

In fact, first of all, notice that that linear cycles of dimension at most $n-1$ spanned by the collection of $s$ points are *Weyl cycles*, according to our definition. This holds because hyperplanes passing through $n$ base points are always Weyl divisors. We recall that for $s = n + 2$, the only Weyl divisors are the exceptional divisors and the hyperplanes spanned by $n$ base points. We conclude that for $s = n+2$ the linear cycles spanned by base points are the only Weyl cycles, hence we have that for divisors in $X_{n+2}^n$, the Weyl expected dimension equals the linear expected dimension, so the analogous of Conjecture (1) in $X_{n+2}^n$ holds by [2]. Moreover, by [2, Corollary 4.8] we can say that the analogous of Conjecture 1 holds in arbitrary dimension for small number of points.

Secondly, in [3] the authors considered cycles $J(L_I, \sigma_t)$, joins over the $t$ secant variety to the rational normal curve of degree $n$ passing through $n+3$ points, and they gave a *secant expected dimension* for an effective divisor. It agrees with the Weyl expected dimension for $n = 4$. For $X_7^4$, these varieties are just the unique rational normal quartic curve through the 7 points and the pointed cones over it, namely cone over rational normal curve, labeled $S^3_{1,\widehat{8}}$ as in notation (1). Therefore, we ask the following question:

**Conjecture 2** *The varieties $J(L_I, \sigma_t)$ are Weyl cycles on $X_{n+3}^n$.*

Solving Conjecture 1 would complete the analysis of the dimensionality problem for all the Mori Dream Spaces of the form $X_s^n$, which are $X_{n+3}^n$, $X_8^2$, $X_7^3$ and $X_8^4$.